\input amstex 
\documentstyle{amsppt} 
\magnification\magstep1
\NoBlackBoxes
\NoPageNumbers
\topmatter  %XXXXXXXXXXXX

\title On the cyclic homology of commutative algebras over arbitrary
ground rings\endtitle
\author by Guillermo Corti\~ nas\endauthor
\leftheadtext{Guillermo Corti\~ nas}
\rightheadtext{Cyclic homology of commutative algebras}
\address Dep. de Matem\'atica, Fac. Cs. Exactas, Calle 50 y 115, 
(1900) La Plata, Argentina\endaddress

\endtopmatter  %XXXXXXXXXXX

\define\am{{*-1}}

\define\dedo{\bar{d}}
\define\defor{\tilde\twoheadrightarrow}

\define\px{{p-1}}

\define\pq{{p,q}}
\define\cris{\bar{\Omega}}

\define\cd{\Cal D}
\define\cf{\Cal F}

\define\ch{\underline{ChCGA}}

\define\de{\delta}

\define\fib{\twoheadrightarrow}

\define\ga{\Gamma}

\define\gp{\gamma^p}

\define\h1{\hat f}

\define\hh{HH_*}

\define\lv{\Lambda V}
\define\lw{\Lambda W}

\define\pa{\partial}
\define\rat{\Bbb Q}

\define\rt{\longrightarrow}

\define\tr{\tilde\rightarrow}
\define\tot{\operatorname{Tot}}

\define\wh{\widehat{HH}_n}
\define\wc{\widehat{HC}_n}

\define\ho{\operatorname{Hom}}

\define\Ker{\operatorname{Ker}}

%XXXXXXXXXXXXXXXXXXXXXXXXXXXXXXXXXXXXXXXXXXXXXXXXX

\subhead{0. Introduction}\endsubhead
\bigskip
We consider commutative algebras and chain algebras over a fixed commutative
ground ring $k$ as in the title. 
We are concerned with the problem of computing the cyclic (and Hochschild)
homology of such algebras via free DG-resolutions $\lv\defor A$. We find
spectral sequences: 
$$
\align
E^2_{p,q}=H_p(\lv\otimes\ga^q(dV))&\Rightarrow HH_{p+q}(\lv)\tag{1}\\
 {E'}^2_{\pq}=H_p(\lv\otimes\ga^{\le q}(dV))&\Rightarrow HC_{p+q}(\lv)
\endalign
$$
The algebra $\lv\otimes\ga(dV)$ is a divided power version of the de 
Rham algebra; in the 
particular case when $k$ is a field of characteristic zero, the spectral
sequences above agree with those found in [BuV], where it is shown they
degenerate at the $E^2$ term. For arbitrary ground rings we prove here 
(theorem 2.3) that if $V_n=0$ for $n\ge 2$ then $E^2=E^\infty$.
From this we derive a formula
for the Hochschild homology of flat complete intersections in terms
of a filtration of the complex for crystalline cohomology, and find a 
description of ${E'}^2$ also in terms of crystalline cohomology (theorem 
3.0).
The latter spectral sequence degenerates for complete intersections of
embedding dimension $\le 2$ (Corollary 3.1).
Without flatness assumptions, our results can be viewed as the computation
Shukla (cyclic) homology (as defined in 1.4 below--see also [PW], [S]). 
Particular cases of theorem 3.0
 and corollary 3.1 have been obtained in [GG],[L] and [LL]. To our knowledge
this is the first paper to give a unified proof for all of these.
\bigskip
The remainder of this paper is organized as follows. Section 1
contains the definitions, results and notations used in the paper. 
Section 2 is devoted to computing the homology of free algebras. The spectral
sequences \thetag{1} are obtained in 2.1-2.2; the degeneracy result is 
proved in 2.3.; the fact that the same degeneracy result is not valid for
general free algebras $\lv$ with  
$V_2\ne 0$ --unless $k\supset\rat$-- is proved in 2.4.
In Section 3 we apply the results of 
section 2 to complete intersections (Theorem 3.0 and Corollary 3.1).

\newpage
\subhead{1. Definitions and notations}\endsubhead
\bigskip
\subsubhead{1.0 Chain algebras}\endsubsubhead
We consider algebras over a fixed {\it ground ring} $k$ which
is assumed unital and commutative. Algebras are unital
and come equipped with a non-negative gradation $A=\oplus_{n\ge 0} A_n$ and a differential $\partial : A_n@>>>A_{n-1}$ with $\partial^2=0$ and $\partial(ab)=\partial(a)b+(-1)^{|a|}a\partial{b}$, 
and are strictly commutative in the graded sense. Thus not only is 
$ba=(-1)^{|a||b|}ab$
but also $a^2=0$ if $|a|$ is odd. Algebras with this added
structure and strict commutativity condition are called
{\it chain algebras} or simply algebras. Usual commutative 
algebras are considered as chain algebras concentrated in degree zero. 
Two algebra
maps $f,g:A@>>>B$ are {\it homotopic} --and we write $f\approx g$--
if they are homotopic
as chain maps, and are {\it specially homotopic} if the
homotopy $s:A_*@>>>B_{*+1}$ is special, i.e. if it verifies $s(ab)=s(a)b+(-1)^{|a|}as(b)$. Algebra maps which induce an
isomorphism in homology are called {\it quisms}
or quasi-isomorphisms and are denoted $@>\sim>>$. Maps which are
surjective in every degree are denoted $\fib$. We abuse language 
and call a chain algebra {\it free} if it is free as a strictly commutative graded algebra. This means
that $A=\lv$ is the symmetric algebra 
--in the graded sense-- of some free graded module $V$; no condition is
imposed on the differential $\partial:A_*@>>>A_{*-1}$. 
\bigskip
The following theorem is well-known
to specialists. The particular case of 1.1-ii) 
when $k$ is a field was proved in [BuV, 1.1].  
\bigskip

\proclaim{Theorem 1.1}(Folklore)

\item{i)} Let $\lv$ be a free algebra,
and let $f\negthickspace:\negthickspace A\defor B$ be a surjective quism. 
Then the induced map
$$
f_*\negthickspace:\negthickspace \ho (\lv,A)\fib\ho (\lv,B)
$$
is surjective.  
Moreover, if $f_*(g^1)\approx f_* (g^2)$ are homotopic
via a special homotopy $s\negthickspace:\negthickspace \lv_*\rt B_{*+1}$, then $g^1\approx g^2$ via
a special homotopy $\hat s\negthickspace:\negthickspace \lv_*\rt A_{*+1}$ such that $f\circ\hat s=s$.
\smallskip
\item{ii)} For every $A\in\ch$ there exists a
free algebra $\lv$ and a surjective quism
$$
\lv\tilde\twoheadrightarrow A
$$
Such a surjective quism will be called a {\it model} of $A$. If
furthermore $A$ has finite type and $k$ is noetherian, then $\lv$ can
be chosen of finite type.
\endproclaim

We omit the proof of the theorem above for lack of space;
it is available upon request.
For the proof of 2.4 below we shall need the following
explicit version of 1.1-ii), which in addition gives an idea of the
proof of the theorem above.
\bigskip 
\proclaim{1.1.1 Addendum} 
Suppose the following data is given:
\smallskip
\item{i)} A free algebra $Q^n=\lv^n$, with
$V^n=V_0\oplus\cdots\oplus V_n$.
\smallskip
\item{ii)} A morphism $f^n\negthickspace:\negthickspace Q^n\rt A$
inducing an isomorphism $H_m(Q^n)@>\cong>>H_m(A)$ for $m<n$
and a surjection
$Z_n(Q^n)=\ker(\partial:Q^n_n@>>>Q^n_{n-1})\fib Z_n(A)$.

Then there exist:

\item{i')} A free chain algebra $Q=\lv$ with $V_i=V^n_i$
for $i=0\dots n$.
\item{ii'')} A surjective quism
$f\negthickspace:\negthickspace Q\tr A$ with
$f\big|_{Q^n}=f^n$.
\endproclaim

\bigskip
\proclaim{Corollary 1.2}
 Let $F\negthickspace:\negthickspace \ch\rt\cd$ be a
functor from the category of chain algebras and with values in any category $\cd$. Assume $F$ maps quisms 
and specially homotopic maps between free algebras to isomorphisms  and to equal maps. Then
the following construction is functorial.
For each $A\in\ch$ choose a model $\epsilon_A:\lv(A)\defor A$ --if 
$A$ is free already, choose $\epsilon_A=id_A$--  and
set $\hat F(A)=F(\lv(A))$. For each map $A@>>>B\in\ch$ choose
a lifting $f':\lv(A)@>>>\lv(B)$ of $f$ with $\epsilon_Bf'=f\epsilon_A$, and set $\hat F(f)=F(f')$. Furthermore the functor $\hat{F}$ maps all
quisms into isomorphisms, the map $F(\epsilon_A):\hat{F}(A)@>>>F(A)$
is a natural transformation, and is universal (final) among all natural
transformations whose source is a functor mapping quisms into isomorphisms.
\endproclaim 
\demo{Proof} Straightforward from Theorem 1.1.\qed\enddemo
\bigskip

\subsubhead{1.3 Double Mixed complexes}\endsubsubhead
By a {\it double mixed complex} we shall understand a chain
complex in the category of (single) mixed complexes in the 
sense of Kassel. Thus a double mixed complex
is a bigraded module $(p,q)\mapsto M_{p,q}$ equipped with three
$k$-linear maps of degree $\pm 1$: $D$ which lowers the $q$ index and fixes $p$;
$\partial$ which fixes $q$ and lowers $p$, and $B$ which increases
$q$ and fixes $p$. These maps satisfy $0=D^2=B^2=\partial^2=D\partial+\partial D=B\partial +\partial B=DB+BD$.
The Hochschild homology of a double mixed complex $M$ as above
is the homology of the double complex $(M_{*,*},D,\partial)$, and its
cyclic homology is the homology of the triple complex $(_BM_{*,*,*},D,\partial, B)$
where for each fixed value of $r$, $_BM_{*,*,r}$ is the usual triangular
double complex for cyclic homology ([K]). The homology of a triple complex
is defined by taking the usual $\tot$ of double complexes twice. One
can do this in different ways which of course yield the same complex
but which suggest different filtrations and thus different spectral sequences. We like to think of $\tot_BM_{*,*}$
as the $\tot$ of the double complex 
$M'_{p,q}=\bigoplus_{i=0}^q M_{p-i,q-i}$ with $D$ as vertical boundary and 
$B+\partial$ as row boundary. Our
choice of spectral sequences is coherent with this. We use the column 
filtration of $M_{*,*}$ to obtain the spectral sequence:
$$
E^1_{\pq}:=H_q(M_{p,*},D)\Rightarrow HH_{p+q}(M)\tag{2}
$$
Note that $E^1=(E^1_{*,*},0,\partial,B)$ is a double mixed complex. We write ${E'}^1$
for the corresponding double
complex as above. Thus ${E'}^1$ is the first term of the spectral sequence
for the column filtration of the double complex $M'$. For the second term in 
the spectral sequence $E$ we 
have to take homology with respect to $\partial$, while for ${E'}^2$ we
take homology with respect to $\partial+B$. Hence
if we put:
$$
HH_n^p(M):=E^\infty_{n-p,p}\text{\quad and\quad} HC_n^p(M):={E'}^\infty_{n-p,p}
\tag{3}
$$
we obtain:
$$
E^2_{p,q}=HH_{p+q}^q(E^1)\Rightarrow HH_{p+q}(M)\text{\quad and\quad }
{E'}^2_{\pq}=HC_{p+q}^q(E^1)\Rightarrow HC_{p+q}(M)\tag{4}
$$
We shall be especially concerned with two double mixed
complexes. One is the cyclic mixed complex $(C_{*,*}(A),b,B,\partial)$
with $C_{p,q}(A)=(A\otimes (A/k)^{\otimes q})_p$, and with the usual 
boundary maps, see e.g [CGG, 1.6]. Another is the
mixed complex of
$\ga$-differential forms, which we shall define in 1.6 below. 
\bigskip
\proclaim{Lemma-Definition 1.4}(Shukla cyclic homology--compare [PW], [S])
 For each $n\ge 0$, the functors $A\mapsto HH_n(A), HC_n(A)$, going from chain
 algebras to filtered modules
equipped with the filtration of 1.3 above, satisfy the hypothesis of Corollary 1.2. In particular
$\wh(A)$ and $\wc(A)$ are defined; we call them 
{\rm Shukla homology} and {\rm Shukla cyclic homology};
the filtration they carry is what we shall call the 
{\rm Hodge filtration}. The $p$-th layers of each filtration
are $\wh^p(A):={HH}_n^p(C(\lv))$ and 
$\wc^p(A):={HC}_n^p(C(\lv))$. 
\endproclaim
\demo{Proof} Immediate from the fact that free algebras are flat as
modules, and the fact that a special homotopy $s$ between maps
$f,g:A@>>>B$ of not necessarily free algebras $A$ and $B$ induces
a homotopy $s :(C_{*,q}(A),\partial)@>>>(C_{*,q}(B),\partial)$ between the maps $f$ and $g$ 
induce on each column of $C_{*,*}$.\qed
\enddemo
\bigskip
\subsubhead{1.5 Divided powers}\endsubsubhead
\bigskip
Recall from [B,I.1.1] that a system of {\it divided power operations}
for a pair $(A,I)$ consisting of a plain algebra
$A$ and an ideal $I\subset A$ is a family of maps
$\gamma_n:I@>>> A$ ($n\ge 0$), with $\gamma_0(x)=1,\quad \gamma_1(x)=x,\quad \gp(x)\in I\quad (p\ge 1)$ and satisfying a number
of identities ([B I.1.1.1-1.1.6]) which make them formally analogue to 
$x\mapsto\frac{x^n}{n!}$.
It is shown in [B, I.2.3.1], that the forgetful functor going from
the category of triples $(A,I,\gamma)$ to the category of pairs $(A,I)$
has a left adjoint. We write $D(A,I)=(\bar{A},\bar{I},\gamma)$ for this adjoint functor, which 
we call the {\it divided power envelope} of $(A,I)$. Oftentimes we
shall abuse notation and write $D(A,I)$ --or simply $D$-- for $\bar{A}$.
The ideals $D\supset
\cf_n:=<\{\gamma^{q_1}(x_1)\dots\gamma^{q_r}(x_r):r\ge 1, 
\sum_{i=1}^rq_i\ge n\}>$ define a descending filtration.
We call $\cf$ the {\it $\gamma$-filtration}. By the {\it crystalline
complex} of $(A,I)$ we mean the largest quotient $\bar{\Omega}$ of the de 
Rham cochain
algebra of $\bar{A}$ for which the induced derivation $\dedo$ is a 
{\it $\gamma$-derivation}, i.e. maps $\gamma_n(x)\mapsto\gamma_{n-1}(x)dx$;
explicitly  $\bar{\Omega}=\Omega/<d\gamma_n(x)-\gamma_{n-1}(x)dx>$. 
There is a natural isomorphism $\bar{\Omega}\cong D\otimes_A\Omega$ of
cochain algebras with $\gamma$-derivation (cf. [I ,Ch.0, 3.1.6]).
Our choice of name for $\bar{\Omega}$ comes from the fact that, 
for example if $k$ is the field with $p$ elements and $A$ is smooth, 
essentially of finite type over $k$, 
it computes the crystalline
 cohomology $H_{crys}^*(A/I)$
of $A/I$ over $k$ (cf. [I,Ch.0, 3.2.3-4.]). 

We shall be especially interested in the particular case of the divided
power envelope of pairs of the form $(\lv,<V>)$ where $V$ is a $k$- module
and $\lv$ is the symmetric algebra. In this case we write $\ga(V)$ for
$D(\lv,<V>)$. If $V$ is a graded $k$-module, we put
$\ga(V):=(\lv_{odd})\otimes\ga(V_{even})$. If $V$ happens to be free on
a homogeneous basis $\{v_i:i\in I\}$, then $\ga(V)$ is free on the
homogeneous basis $v_I\gamma^Q(v_J):=v_{i_1}\dots
v_{i_r}\gamma^{q_1}(v_{j_1})\dots\gamma^{q_s}(v_{j_s})$, where
$I=(i_1,\dots i_r)$, $J=(j_1,\dots j_s)$, $Q=(q_1\dots q_s)$, the $v_i$
are of odd degree and the $v_j$ are of even degree. We write $\ga^q(V)$
for the $k$-submodule generated by all the $v_I\gamma^Q(v_J)$ with
$q=|Q|:=\sum_{i=0}^sq_i$, and $\ga^{\le q}(V)=\oplus_{^p=0}^q\ga^p(V)$.
\bigskip 
\subsubhead{1.6. The mixed complex of $\ga$-forms}\endsubsubhead
Let $\lv$ be a free chain algebra, and let $\partial$ be its differential.
Write $dV$ for the graded module $V$ shifted by one: $dV_n=V_{n-1}$. Form
the graded algebra $\lv\otimes\ga(dV)$ and equip it firstly with the
degree-increasing $\gamma$-derivation $d$ extending $v\mapsto dv$ and
secondly with the degree-decreasing $\gamma$-derivation $\delta$ with 
prescriptions $\delta(a)=\partial(a)$ ($a\in \lv$) and  
$\delta{dv}=-d\partial(v)$. The (double) {\it mixed complex of 
$\ga$-differential forms} is
$M(\lv)_{p,q}:=(\lv\otimes\ga^q(V))_p$ with $0$ as the $D$ boundary,
$\delta$ as the $\partial$ boundary and $d$ as the $B$ boundary. Note that
the spectral sequences \thetag{4} for this complex degenerate at the
second term. We shall show in 2.1 below that $M(\lv)$ is isomorphic to the
spectral mixed complex $E^1$ of $C(\lv)$. 
\bigskip
\subhead{2. The homology of free algebras}\endsubhead
\bigskip
\proclaim{Lemma 2.0} Let $\lv=(\lv,0),\  \lw=(\lw,0)$ be free
algebras, with zero boundary maps.  For each $p\geq 0$, let
$K^p=\lv\otimes\lw\otimes\ga^p(dW)$, and put:
$$
K=\oplus_{p\geq 0}K^p\simeq\lv\otimes\lw\otimes\ga(dW)
$$
Equip $K$ with the $\gamma$-derivation $D\negthickspace:\negthickspace K_*^p\rt K_*^{\px}$
determined by $D(dw)=w$

\noindent  $(w\in W)$, $D^2=0$.
Then:
\roster
\item"{i)}" Consider the natural projection
$\pi:\lv\otimes\lw@>>>\lv$; then $\ker\pi=<W>=D(K^1)$ is both the
ideal generated by $W$ and the image of $D$. We write $I$ for this ideal. 
\smallskip
\item"{ii)}" The ideal $K'=\Ker(\pi\otimes Id)\subset K$ is
contractible.  Precisely, there is a $\lv$-linear map
$h\negthickspace:\negthickspace K_*^{\prime p}\rt K_{*+1}^{\prime p +1}$
with
$hD+Dh=Id,\ h(w)=dw,$ $ h^2=0$ $(w\in W)$
and 
$h(I^r K)\subset I^{r-1} K$ $(r\geq 1)$. In particular
$\pi\otimes 1\negthickspace:\negthickspace K\tilde\twoheadrightarrow\lv$
is a free resolution of $\lv$ as a $\lv\otimes\lw$-module.
\endroster
\endproclaim

\smallskip
\demo{Proof} Part i) is trivial.  Next observe that it suffices to prove ii) 
in 
the case $V=0$, for if 
$h=h(0):\lw\otimes\ga^*(dW)\mapsto\lw\otimes\ga^{*+1}(dW)$ is a homotopy
satisfying the prescriptions of the lemma for $V=0$, then 
$h(V)=1\otimes h:\lv\otimes\lw\otimes\ga^*(dW)@>>>\lv\otimes\lw\otimes\ga^{*+1}
(dW)$ satisfies the prescriptions for $V$. Assume $V=0$; if $W=0$ there
is nothing to prove. Otherwise choose
a well-ordered, nonempty basis $(B,<)$ of $W$. For each $w\in B$, define
$h:\Lambda(kw)\otimes\Gamma^*(kdw)@>>>\Lambda(kw)\otimes\Gamma^{*+1}(kdw)$
as follows.  Put $h(w^ndw)=0$ and  $h(w^{n+1})=w^ndw$ if $|w|$ is even,
$h(w\gamma_p(dw))=\gamma_{p+1}(dw)$ and $h(\gamma_p(dw))=0$ if $|w|$ is
odd, and $h(1)=0$ in either case. Given any strictly increasing sequence
$w_1<\dots<w_r$, extend $h$ to 
$A(w_1,\dots, w_r):=\Lambda(kw_1)\otimes\Gamma^{*+1}(kdw_1)\otimes\dots\otimes 
\Lambda(kw_r)\otimes\Gamma^{*+1}(kdw_r)$ by $h(x_1\otimes\dots\otimes x_r)=
(-1)^{|x_1|+\dots |x_{r-1}|}x_1\otimes\dots x_{r-1}\otimes h(x_r)$. One checks
immediately that these prescriptions give a well defined linear map
$h: \lw\otimes\ga^*(\lw)=\bigcup_{r\ge 1}
\bigcup_{w_1<\dots<w_r\in B}A(w_1,\dots, w_r)@>>>\lw\otimes\ga^{*+1}(\lw)$
and that this map is a homotopy satisfying the requirements of the lemma.\qed
\enddemo
\bigskip
\proclaim{Proposition 2.1} Let $\lv$ be a free chain algebra. Consider
the spectral sequence $E^1$ of \thetag{2} associated to the standard mixed 
complex $C(\lv)$. Then the double mixed complex $(E^1_{*,*}, 0,\partial, B)$
is isomorphic to the double mixed complex of $\ga$-differential forms defined 
in 1.6 above. 
\endproclaim 
\demo{Proof} The plan of the proof is as follows. Firstly we construct a 
map $f:\lv\otimes\ga(dV)@>>>E^1$ and show it commutes with the relevant
boundary maps. Secondly we find an isomorphism 
$\bar F:\lv\otimes\ga(dV)\cong E^1$. Thirdly we prove that $f=\bar F$.
For the definition of $f$ proceed as follows. It is not hard to see
that $d:\lv @>>>\lv\otimes\ga^1(dV)$ is the universal derivation, so that
$\lv\otimes\ga^1(dV)$ is the graded $\lv$-module of K\"ahler differentials
of $\lv$. Hence $\lv\otimes\ga^1(dV)\cong E^1_{*,1}$ as graded modules,
and the isomorphism maps $da$ to the class of $Ba$ ($a\in\lv$). Hence 
$\partial(da)=-d\partial(a)$ goes to the class of $-B(\partial(a))=\partial(B(a))$,
whence the isomorphism $\lv\otimes\ga^1(dV)\cong E^1_{*,1}$ is a chain
module map, because it is a module map which commutes with the boundary on
the generators $da$. Next, the shuffle product and divided power operations of Cartan
([Ca, Exp.4, Th.4; Exp.7, Th.1]) --actually their analogue for the cyclic
bar construction-- make $E^1$ into a DG algebra with divided operations for
the ideal generated by $E^1_{*,1}$. By the universal property of
$\lv\otimes\ga(dV)$, we have a map $f:\lv\otimes\ga(dV)@>>>E^1$ commuting
with $\partial$. Next we must check that $fd=Bf$. Note that the image of 
$f$ is generated, as a graded $k$-module, by classes of monomials of the 
form $M=f(a\gamma_{p_1}(dy_1)\dots\gamma_{p_r}(dy_r))=
[a*(1\otimes y_1^{\otimes p_1})*\dots *(1\otimes y_r^{\otimes p_r})]$, where 
$[$ $]$ denotes
homology class, $*$ is the shuffle
product, $a\in\lv$, $y_i\in V$, $r\ge 1$, $p_i\ge 0$ and $p_i\le 1$ if $|y_i|$ is
even. Thus we must show that 
$B(M)=[Ba*(1\otimes y_1^{\otimes p_1})*\dots *(1\otimes y_r^{\otimes p_r})]$.
The same proof as in the ungraded case [LQ, Lemma 3.1] proves that the formula
$B(x*By)=Bx*By$ holds for chain algebras. Using this
formula and induction, we see it suffices to show 
that, for $y\in \lv_{even}$, $x\in C(\lv)$ and $p\ge 2$, we have 
$$
B(x*(1\otimes y^{\otimes p}))=B(x)*(1\otimes y^{\otimes p})\tag{5}
$$
Whenever $p$ is not a zero divisor in $k$, \thetag{5} follows from the
formula of [LQ, 3.1] and the fact that $1\otimes y^{\otimes p}=\frac{1}{p}B(y^{\otimes p})$
over $k[1/p]$. In particular \thetag{5} holds for $\Bbb Z$, whence it
holds for arbitrary ground rings, by naturality. This finishes the
first part of the proof. Next apply the lemma above with
$1\otimes V$ as $V$ and $TV=\{ v\otimes 1-1\otimes v:v\in V\}$ as $W$,
 to obtain a free resolution
$(K,0)=(\lv\otimes \lv\otimes\ga(dV),D)\tilde\twoheadrightarrow (\lv, 0)$.
Write $C'$ for the Hochschild acyclic resolution.
By universal property, the inclusion map
$F\negthickspace:\negthickspace\lv\otimes\lv\otimes\ga^1(dT V)\hookrightarrow
C^{\prime 1}(\lv),\ \
dTv\longmapsto 1\otimes v \otimes 1$ extends to an algebra 
homomorphism $F\negthickspace:\negthickspace K\rt C'(\lv)$
with $F(K^p)\subset C^{\prime p}(\lv)$ and $F(\gamma^p(dTv))=\gp F(dTv)$, 
$p\geq 0$, $v\in V$. Because $F(DdTv)=Tv=b'F(dTv)$, $F$ is a chain map.  
Since both $K$ and $C'$ are resolutions, the induced map
$\bar F\negthickspace:\negthickspace \bar K=(\lv\otimes\ga(dV), 0)\tilde{\rt} C(\lv), b)$
is a quism.  Thus, $\bar F$ induces an isomorphism of divided power
algebras $\ga_q^p(\lv)\simeq E_\pq^1$.
But $\bar F(dv)=f(dv)$ for all $v\in V$. Thus, $\bar F$
and $f$ are the same map.$\qed$\enddemo
\bigskip
\proclaim{Corollary 2.2} Let $A$ be an algebra, $\lv\defor A$
be a free model, $M(\lv)$ the mixed complex of $\gamma$-forms of 1.6
above. Then 
there are spectral sequences:
$$
E^2_{\pq}=HH_{p+q}^p(M(\lv))\Rightarrow \widehat{HH}_{p+q}(A),\quad 
{E^\prime}^2_{p,q}=HC^p_{p+q}(M(\lv))\Rightarrow \widehat{HC}_{p+q}(A)
$$
Furthermore, for the Hodge filtration, we have 
$E^\infty_{p,q}=\widehat{HH}_{p+q}^p(A)$ and 
${E^\prime}^\infty_{p,q}=\widehat{HC}_{p+q}^p(A)$.
\endproclaim
\demo{Proof} Immediate from lemma 1.4.\qed
\enddemo
\bigskip
\proclaim{Theorem 2.3} In the situation of the proposition above, assume
further that $V_n=0$ for $n\ge 2$.  Then the spectral sequence for 
Hochschild homology degenerates, and there is an isomorphism {\it of 
graded algebras}
$$
\hh(\lv)\cong H_*(\lv\otimes\ga(dV))
$$
\endproclaim
\demo{Proof} We use the notations of the proof of the proposition above.
The plan of the proof is as follows.  Firstly, we equip the algebra $K$
with a $\gamma$-derivation $\Delta\negthickspace:\negthickspace K\rt K$. 
Secondly, we define a chain homomorphism 
$F'\negthickspace:\negthickspace (K,\Delta)\rt (C', b'+\pa')$.  
Thirdly, we show that the induced map 
$\bar F'\negthickspace:\negthickspace (\bar K,\bar \Delta)\rt(C,b +\pa)$ 
is a quism.  Finally, we prove that $(\bar K,\bar \Delta)\simeq
(\lv\otimes\ga(dV),\de)$.  Let $\de'$ be the $\gamma$-derivation determined
by
$\de'(x)=\pa'x$, $\de'(dTv_0)=0$, $\de'(dTv_1)=-h(\pa v_1)$, $x\in
\lv\otimes \lv$, $v_i\in V_i$, $i=0,1$. By definition of $\de'$, we have
$\de^{\prime 2}=0$ and $(D+\de')^2=0$.
Thus $\Delta:=D+\de'\negthickspace:\negthickspace  K_*\rt K_{\am}$
is a boundary derivation, and $(K,\Delta)$ is a chain algebra. 
Let $F\negthickspace:\negthickspace K\rt C'$ be as in the proof
of the proposition above; then $F$ comutes with
$D$ and $b'$ but not necessarily with $\Delta$. Define a morphism
$F'\negthickspace:\negthickspace K\rt C'$
as follows. Firstly, let
$F'(w)=F(w)$, $w\in L=\lv\otimes \lv\otimes\ga(dTV_0)\subset K$.
Secondly, choose a basis $\Cal B$ of $V_1$ and put
$F'(dTv)=F(dTv) +1\otimes ((F\de'- \partial' F)dT(v)))$ ($v\in \Cal B$).
Because $x\longmapsto 1\otimes x$ is a chain homotopy, we have:
$$
F'(\Delta dTv)=(b'+\pa') F'(dTv)\tag{6}
$$
We remark that 
$G(dTv)=(F'- F)(dTv)\in \lv_0\otimes \lv_0\otimes \ga^2(TV_0)$.
In particular:
$$
(G(dTv))^2=0\qquad (v\in \Cal B)\tag{7}
$$
For each $p\geq 1$, put
$F'(\gamma^p(dTv))=\gamma^p F(dTv)+\gamma^\px F(dTv) G(dTv)$ ($v\in\Cal B$).
Fix $v\in\Cal B$, and write $\eta=dTv$. 
Then:
$$
\align
F'(\gp(\eta)) F'(\gamma^q(\eta)) &= \\
{\text{(by \thetag{7})}} &=\gp F(\eta)\gamma^q F(\eta) +
(\gamma^{\px} F(\eta) \gamma^q F(\eta) + \gp F(\eta)\gamma^{q-1}
F(\eta)) G(\eta) \\
 &=\binom{p+q}p \gamma^{p+q} F(\eta) +
\biggl[\binom{p+q-1}{p-1} + \\ 
&\quad\binom{p+q-1}p\biggr]\gamma^{p+q-1}(F(\eta))G(\eta)\\
 &=\binom{p+q}p F(\gamma^{p+q}(\eta))
\endalign
$$
Thus, the maps $\eta\longmapsto F(\gp(\eta)),\ p\geq 1$, extend to
divided power operations
$\theta^p\negthickspace:\negthickspace dT(V)\rt C',\ \ v\longmapsto F(\gp(dTv))$
($v\in \Cal B$)
satisfying all of the conditions of [B I.1.1].  By \thetag{6}, $b' +\pa'$ 
is a $\theta$-derivation.  Thus, $F'$ extends uniquely to a chain
homomorphism
$F'\negthickspace:\negthickspace K\rt C'$ with $F^\prime \circ \gp(\eta) 
=\theta^p\circ F'(\eta)$, $p\geq 1$. Let
$\bar F\negthickspace:\negthickspace (K,\Delta)\otimes_{\lv\otimes\lv} \lv\simeq
(\lv\otimes\ga(dV),\bar\de')\rt (C, b+\pa)$
be the induced map.  Consider the subalgebra $L\subset K$ defined
above.  By definition, $\bar F'=\bar F$ both on $L$ and on $K/L$. 
Since $\bar F$ is a quism, so is $\bar F'$.  To finish the proof, it
suffices to show that
$\bar \de'=\de\negthickspace:\negthickspace \lv\otimes\ga(dV)_*\rt \lv\otimes\ga(dV)_\am$
is the usual derivation.  Let $\Cal B_*\subset V_*$ be the ordered basis
used in the definition of the homotopy $h\negthickspace:\negthickspace K\rt K$.  
Because both
$\bar\de'$ and $\de$ are $\gamma$-derivations and $\bar\de'\equiv\de$ on
$L\subset K$, it is enough to show that:
$$
\bar\de' (dv)= -\mu(hT(\pa v))=-d\pa v=\de(dv)\qquad (v\in
\Cal B_1)\tag{8}
$$
Here $\mu\negthickspace:\negthickspace C'\twoheadrightarrow C$ is the projection map
and $T\negthickspace:\negthickspace \lv\rt
J=\ker(b':\lv\otimes\lv\rt\lv)$,
$T(a)=a\otimes 1-1\otimes a$
is the universal non commutative derivation of degree zero.
We have $\Ker \mu=J C'$ and:
$$
T(ab)\equiv(a\otimes 1)T(b)+(b\otimes 1)T(a)\qquad({\text{mod}}\  J^2)\tag{9}
$$
Now, for each $v\in B_1$, $\pa v$ is a sum of monomials of the form:
$\lambda v_1\ldots v_r,\ \ \lambda\in k,\ \ v_1<\ldots< v_r,\ \ v_i\in
V_0$.
Thus $d\pa v$ is the sum of the monomials
$\lambda v_1\ldots\overset{v_i}\to \vee\ldots v_r dv_i$.
On the other hand, by \thetag{9}, $T\pa v$ is the sum of the monomials
$\lambda v_1\ldots\overset{v_i}\to \vee\ldots v_r Tv_i$
plus an element in $J^2$.  But by Lemma 2.0-ii), we have $h(J^2)\subset J 
K^1$.
Thus $h(T\pa V)$ is congruent modulo $JK^1$ to the sum of the monomials
$\lambda v_1\ldots\overset{v_i}\to \vee\ldots v_r dTv_i$.
It follows that $\mu (T\pa v)=d(\pa v)$.
We have established the identity \thetag{8}; this concludes the
proof.$\qed$\enddemo
\bigskip
\proclaim{Proposition 2.4} Let $k$ be any ground
ring.  Suppose that for every model
$\lv \tilde\twoheadrightarrow k$ ($V$ arbitrary)
the spectral sequence $E$ of Corollary 2.2 degenerates at the second
term. Then $k\supset\Bbb Q$.
\endproclaim
\demo{Proof} Suppose that $k\not\supset\Bbb Q$; then there is an integer 
$p\geq 2$ that is not invertible in $k$.  We are going to exhibit a
free algebra $\lv$ that is quasi-isomorphic to $k=\Lambda(0)$, but for 
which $H_{2p}(\lv\otimes\ga(dV))\neq 0=H_{2p}(\Lambda(0)\otimes\ga(0))$.
Consider the graded $k$- module $V'$ defined by $V'_0=V'_n=0$ if $n\ne 1,2$
and $V'_1=ky$, $V'_2=kz$. Equip $\lv'$ with the derivation $\partial'y=0$,
$\partial'z=y$. Let $f':\lv'@>>>k$ be the natural projection. 
Thus, $H_0(\lv')=k$ and $H_i(\lv')=0$ for $i=1,2$.
Hence by 1.1.1, $f'$ extends to a quism $f:\lv\defor k$ where 
$V_i=V'_i$ for $i\le 2$ and where $\partial(a')=\partial'(a')$ 
for $a'\in \lv'$. Consider the element 
$\gamma^p(dy)\in(\lv\otimes\ga^p(dV))_{2p}$; 
we have $\delta(\gp(dy))=\gamma^\px(dy)(-d\pa y)=0$.
I claim that $\gp(dy)$ is not a boundary.  Suppose otherwise that
there exists $\alpha\in(\lv\otimes\ga^p(dV))_{2p+1}$ with 
$\de\alpha=\gp(dy)$. 
Extend $\{y,z\}$ to a basis of $V$. Then $(\lv\otimes\ga^p(dV))_{2p+1}$
is the free
$\lv$-module on the monomials $v_I\gamma^Q(v_J)$ of 1.6 above. 
Because $V_2=kz$, it follows that $\beta=\gamma^{\px}(dy)dz$ is the
only basis element of $\ga_{2p+1}^p(\lv)$ whose boundary is a multiple
of $\gp(dy)$.  Thus $\de(\alpha)$ must be a multiple of $\de(\beta)$. 
But $\de(\beta)=\de(\gamma^\px(dy)dz)=-\gamma^\px(dy)dy=-p\gp(dy)$;
since $p$ is not invertible in $k$, it follows that the element
$\alpha$ does not exist.\qed\enddemo
\bigskip
\subhead{3. Complete Intersections}\endsubhead
\bigskip
\proclaim{Theorem 3.0}({\it Complete Intersections})
Let $I$ be an ideal in the polynomial ring $R=k[x_1,\dots x_n]$,
and set $A=R/I$. 
Assume $I$ is generated by an $R$-sequence, i.e. assume there is
a model $\lv\defor A$ with $V_0=<x_1,\dots , x_n>$ and $V_n=0$ for
$n\ge 2$. Then, with the notations of 1.5 above:
\roster
\item Consider the complex:
$$
L^p:\quad\cris^p/\Cal F_1\cris^p@<\dedo<<\Cal F_1\cris^{p-1}/\Cal
F_2\cris^{p-1}@<\dedo<<\dots@<\dedo<<\Cal F_{p}\cris^0/\Cal F_{p+1}
$$
Then $\wh^p(A)=H_{n-p}(L^p)$ and $\wh(A)\cong\oplus_{p=0}^n\wh^p(A)$.
\medskip
\item Consider the complex:
$$
{L'}^p:\quad \cris^p/\Cal F_1\cris^p@<\dedo<<\cris^{p-1}/\Cal 
F_2\cris^{p-1}@<\dedo<<\dots@<\dedo<<\cris^0/\Cal F_{p+1} 
$$
Then the spectral sequence of Corollary 2.2 has ${E'}^2_{\pq}=H_q({L'}^p)$. 
\endroster
\endproclaim
\demo{Proof} In view of corollary 2.2 and Theorem 2.3, it suffices to
show that there are quasi-isomorphisms $(\lv\otimes\ga^p(V),\delta)\defor
L^p$ and
$(\lv\otimes\ga^{\le p}(dV),d+\delta)\defor {L'}^p$. By virtue of [B,
3.4.4 and 
3.4.9], the argument given in [CGG, proof of Th. 3.3 on page 229] to 
prove the case $k\supset\rat$ works here 
also.\qed\enddemo
\bigskip
\proclaim{Corollary 3.1} If in the theorem above the number of variables
is $r\le 2$ then 
$$
\wc^p(A)=H_{n-p}({L'}^p)\qquad\text{for all\quad} n\ge p\ge 0\qed
$$\endproclaim
\medskip
\remark{Remark 3.2} It is not hard to see that our Hodgewise complexes
$L^p$
of a complete intersection are quasi-isomorphic to 
those found in [GG]. For $r=1,2$ the sum of the Hodgewise complexes
${L'}^p$ of the corollary above give the complexes of [LL, 2.5.3] and
[L, 2.10]. All of these are generalizations of the Feigin-Tsygan complexes
([FT], [CGG]).
\endremark
\bigskip
\remark{Acknowledgements} I wrote a first version of this paper in 1993;
I was just learning the subject of differential graded algebras and
thus felt the need to write down every single detail. The result was
a 75 page manuscript. The project of re-writing this paper into a
publishable form has been first on my procrastination list ever since 
then. I would never have 
done it but for the encouragement and useful suggestions I received from
several people, including L. D'Alfonso, T. Pirashvili, M. Vigu\'e and
C. Weibel; thanks to all of them. Thanks also to the people at the 
Mathematics deparment of the Universidad Nacional del Sur, Bah\'\i a Blanca,
Argentina, for their hospitality. It was during a weeklong visit to that
institution that I first managed to trim down the original manuscript 
to an 11 page draft. The referee's suggestions for 
shortening the proof of Lemma 2.0 helped trim the paper further 
to its current size. 
\endremark


\bigskip

\Refs
\widestnumber\key{WWW}

\ref\key B\by P. Berthelot\book Cohomologie cristalline
des Sch\'emas de caract\'eristique $p>0$.\bookinfo Lecture
Notes in Math\vol 407\publ Springer\publaddr Berlin
\yr 1974\endref

\ref\key BuV\by D. Burghelea and M. Vigu\' e\book Cyclic homology of
commutative algebras I\bookinfo Lecture Notes in Math. {\bf 1318}
\publ Springer Verlag\publaddr Berlin\yr 1988\endref

\ref\key Ca\by H. Cartan\paper Constructions multiplicatives (Expos\'e 
4); Puissances divis\'ees (Expos\'e 7) \inbook Alg\` ebres 
d'Eilenberg-MacLane et
homotopie\publ Seminaire Cartan\vol 7$^e$ ann\' ee \yr 1954-55\endref

\ref\key CGG\by G. Corti\~ nas, J.A. Guccione and J.J.
Guccione\paper Decomposition of the Hochschild and cyclic homology of
commutative differential graded algebras\jour J. Pure Appl. Algebra\vol 83\yr
1992\pages 219--235\endref

\ref\key FT\by B. Feigin and B. Tsygan\paper Additive $K$-theory and
crystalline homology\jour Functional Anal. Appl.\vol 19\yr
1985\pages 124--215\endref

\ref\key GG\by J.A. Guccione and J.J. Guccione\paper Hochschild
homology of complete intersections\jour J. Pure Appl. Algebra\vol 74\yr
1991\pages 159--176\endref

\ref\key I\by L. Illusie\paper Complexe de de Rham-Witt
et cohomologie cristalline\jour Ann. Sci. \'Ecole 
Normale Sup\'erieure\vol 12\yr 1979\pages 501-661\endref

\ref\key K\by C. Kassel\paper Cyclic homology, comodules and mixed
complexes\jour J. of Alg.\vol 107\yr 1987\pages 195--216\endref

\ref\key L\by M. Larsen\paper Filtrations, mixed complexes, and cyclic
homology in mixed characteristic\jour $K$-theory\vol 9\pages 173-198
\yr 1995\endref

\ref\key LL\by M. Larsen, A. Lindenstrauss\paper Cyclic homology
of Dedekind domains\jour $K$-theory\vol 6\yr 1992\pages 301--334\endref

\ref\key LQ\by J.L. Loday, D. Quillen\paper Cyclic homology and the Lie 
algebra homology of matrices\jour Comm. Math. Helv.\vol 59\yr 1984
\pages 565-591\endref

\ref\key PW\by T. Pirashvili, F. Waldhausen\paper Mac Lane homology and
topological Hochschild homology\jour 
J. Pure Appl. Algebra\vol 82\yr 1992\pages 81-98\endref

\ref\key S\by V. Shukla\paper Cohomologie des Alg\` ebres
Associatives\jour Ann. Sci. Ec. Norm. Sup. 3$^e$ s\` erie, t\vol 78\yr
1961\pages 163--209\endref

\endRefs
\enddocument